\documentclass{article}
\usepackage{color}
\usepackage{amsmath,amssymb}
\catcode`\@=11 \@addtoreset{equation}{section}

\catcode`\@=12

\newtheorem{Theorem}{Theorem}[section]
\newtheorem{Proposition}{Proposition}[section]
\newtheorem{Lemma}{Lemma}[section]
\newtheorem{Corollary}{Corollary}[section]
\newtheorem{Definition}{Definition}[section]

\newcommand{\bTheorem}[1]{
\begin{Theorem} \label{T#1} }
\newcommand{\eT}{\end{Theorem}}

\newcommand{\bProposition}[1]{
\begin{Proposition} \label{P#1}}
\newcommand{\eP}{\end{Proposition}}

\newcommand{\bLemma}[1]{
\begin{Lemma} \label{L#1} }
\newcommand{\eL}{\end{Lemma}}

\newcommand{\bCorollary}[1]{
\begin{Corollary} \label{C#1} }
\newcommand{\eC}{\end{Corollary}}

\newcommand{\bFormula}[1]{
\begin{equation} \label{#1}}
\newcommand{\eF}{\end{equation}}

\newcommand{\Ov}[1]{\overline{#1}}
\newcommand{\DC}{C^\infty_c}

\newcommand{\vr}{\varrho}
\newcommand{\vre}{\vr_\ep}
\newcommand{\vte}{\vt_\ep}
\newcommand{\vue}{\vu_\ep}

\newcommand{\vt}{\vartheta}
\newcommand{\vu}{\vc{u}}
\newcommand{\vc}[1]{{\bf #1}}

\newcommand{\Div}{{\rm div}_x}
\newcommand{\Grad}{\nabla_x}

\newcommand{\tn}[1]{\mbox {\F #1}}
\newcommand{\dx}{{\rm d} {x}}
\newcommand{\vph}{{\boldsymbol \varphi}}
\newcommand{\vq}{\mathbf{q}}
\newcommand{\dt}{{\rm d} t }

\newcommand{\dxdt}{\dx\dt}

\newcommand{\K}{{\cal K}}
\newcommand{\ep}{\varepsilon}

\font\F=msbm10 scaled 1000

\providecommand{\abs}[1]{\left\lvert#1\right\rvert}
\newcommand{\sil}{\rightarrow}

\definecolor{grey}{rgb}{0.85,0.85,0.85}
\date{}

\long\def\greybox#1{%
    \newbox\contentbox%
    \newbox\bkgdbox%
    \setbox\contentbox\hbox to \hsize{%
        \vtop{
            \kern\columnsep
            \hbox to \hsize{%
                \kern\columnsep%
                \advance\hsize by -2\columnsep%
                \setlength{\textwidth}{\hsize}%
                \vbox{
                    \parskip=\baselineskip
                    \parindent=0bp
                    #1
                }%
                \kern\columnsep%
            }%
            \kern\columnsep%
        }%
    }%
    \setbox\bkgdbox\vbox{
        \color{grey}
        \hrule width  \wd\contentbox %
               height \ht\contentbox %
               depth  \dp\contentbox
        \color{black}
    }%
    \wd\bkgdbox=0bp%
    \vbox{\hbox to \hsize{\box\bkgdbox\box\contentbox}}%
    \vskip\baselineskip%
}

\begin{document}


\title{Weak solutions to the full Navier-Stokes-Fourier system with slip boundary conditions in time dependent domains}

\author{
Ond{\v r}ej Kreml $^{1}$\thanks{The work of O.K., V.M. and \v S.N. was supported by Grant of GA \v CR GA13-00522S and by RVO 67985840.}
\and V{\' a}clav M{\' a}cha $^1$\footnotemark[1]
\and {\v S}{\' a}rka Ne{\v c}asov{\' a} $^1$\footnotemark[1]
\and Aneta Wr\'oblewska-Kami\'nska $^2$\thanks{The work of A.W.-K. was supported by Grant of National Science Center Sonata, No 2013/09/D/ST1/03692.}
}

\maketitle

\bigskip

\centerline{$^1$Institute of Mathematics of the Academy of Sciences of
the Czech Republic} \centerline{\v Zitn\' a 25, 115 67 Praha 1,
Czech Republic}
\bigskip

\centerline{$^2$  Institute of Mathematics, Polish Academy of Sciences}
\centerline{\'Sniadeckich 8, 00-656 Warszawa, Poland}

\medskip

\begin{abstract}

We consider the compressible Navier-Stokes-Fourier system on
time-dependent domains with prescribed motion of the boundary,
supplemented with
slip boundary conditions for the velocity. Assuming that the pressure
can be decomposed into an elastic part and a thermal part, we prove
global-in-time existence of weak solutions. Our approach is based on the
penalization of the boundary behavior, viscosity, and the pressure in
the weak formulation. Moreover, the thermal energy equation is in the
weak formulation replaced by the thermal energy inequality
complemented with the global total energy inequality. In the approximation
scheme the thermal energy inequality is consider to be satisfied in the renormalized sense.

\end{abstract}

\medskip

{\bf Keywords:} compressible Navier-Stokes-Fourier equations,
time-varying domain, slip boundary conditions

\medskip

\section{Introduction}\label{i}

The flow of a compressible viscous heat conducting fluid is in the absence of external forces described by the following system of partial differential equations
\bFormula{i1a}
\partial_t \vr + \Div (\vr \vu) = 0,
\eF
\bFormula{i1b}
\partial_t (\vr \vu) + \Div (\vr \vu \otimes \vu) + \Grad p(\vr,\vt) = \Div \tn{S}(\Grad \vu),
\eF
\bFormula{i1c}
\partial_t (\vr E) + \Div ((\vr E + p) \vu)+\Div \vq = \Div(\tn{S}(\Grad\vu)\vu).
\eF
These equations are mathematical formulations of the balance of mass, linear momentum and total energy respectively. Here $\vr$ is the density of the fluid, $\vu$ denotes the velocity and $E = \frac12 \abs{\vu}^2 + e(\vr,\vt)$ the specific total energy given as sum of the kinetic energy and the internal energy $e$ being a function of density and temperature $\vt$.

The stress tensor $\tn{S}$ is determined by the standard Newton rheological law
\bFormula{i4}
\tn{S} (\Grad \vu) = \mu \left( \Grad \vu +
\Grad^t \vu - \frac{2}{3} \Div \vu \tn{I} \right) + \eta \Div \vu
\tn{I},\,\, \mu > 0,\ \eta \geq 0.
\eF
For simplicity, throughout the rest of this paper we assume the viscosity coefficients $\mu$ and $\eta$ to be constant, however the case of $\mu, \eta$ being dependent on temperature in a suitable way can be also treated, see Section \ref{d}.

The Fourier law for the heat flux $\vq$ has the following form:
 \bFormula {q1}
 \vq= -\kappa (\vt)\Grad\vt, \,\,
\kappa
>0.
\eF

Motivated by \cite{EF70} we assume the following state equation for the pressure
\bFormula{p1p}
p(\vr,\vt) = p_e(\vr)+\vt p_{\vt}(\vr),
\eF
with the additional assumption that $p_{\vt}$ is a non-decreasing function of the density vanishing for $\vr =0$. Consequently the Maxwell relation yields the form of the specific internal energy as
\bFormula{p1pa}
e(\vr,\vt) = P_e(\vr)+Q(\vt)
\eF
with the elastic potential
\bFormula{p1pb}
P_e(\vr) = \int_1^\vr \frac{p_e(z)}{z^2} {\rm d} z
\eF
and the thermal energy being related to the specific heat at constant volume $c_v$ by
\bFormula{p1pc}
Q(\vt) = \int_0^\vt c_v(z) {\rm d} z, \,\, c_v(z) \geq \underline{c_v}>0 \mbox { for all } z\geq
0,
\eF

Assuming smoothness of the flow the system \eqref{i1a}-\eqref{i1c} can be rewritten using the equation for the thermal energy instead of the balance of total energy
\bFormula{i2}
\partial_t \vr + \Div (\vr \vu) = 0,
\eF
\bFormula{i3}
\partial_t (\vr \vu) + \Div (\vr \vu \otimes \vu) + \Grad p(\vr,\vt) = \Div \tn{S}(\Grad \vu),
\eF
\bFormula{i3a}
\partial_t (\vr Q(\vt)) + \Div (\vr Q(\vt) \vu)+\Div \vq = \tn{S}(\Grad\vu): \Grad \vu - \vt p_{\vt}(\vr)\Div \vu.
\eF

We study the system of equations \eqref{i2}-\eqref{i3a} on a moving domain $\Omega = \Omega_t$ with the prescribed movement of the boundary and on time interval $[0,T]$ with $T<\infty$. More precisely, the boundary of the domain $\Omega_t$ occupied by the fluid is described by means of a \emph{given} velocity field $\vc{V}(t,x)$, where $t \geq 0$ and
$x \in \mathbb{R}^3$. Assuming $\vc{V}$ is regular, we solve the associated system of differential equations
\begin{equation}
\frac{{\rm d}}{{\rm d}t} \vc{X}(t, x) = \vc{V} \Big( t, \vc{X}(t, x) \Big),\ t > 0,\ \vc{X}(0, x) = x,
\end{equation}
and set
\[
\Omega_\tau = \vc{X} \left( \tau, \Omega_0 \right), \ \mbox{where} \ \Omega_0 \subset \mathbb{R}^3 \ \mbox{is a given domain},\
\Gamma_\tau = \partial \Omega_\tau, \ \mbox{and}\ Q_\tau = \{ (t,x) \ |\ t \in (0,\tau), \ x \in \Omega_\tau \}.
\]

The impermeability of the boundary of the physical domain is described by the condition
\bFormula{i6} (\vu - \vc{V})
\cdot \vc{n} |_{\Gamma_\tau} = 0 \ \mbox{for any}\ \tau \geq 0,
\eF
where $\vc{n}(t,x)$ denotes the unit outer normal vector to the boundary $\Gamma_t$.
Moreover, we assume the  Navier type 
boundary conditions in
the form
\bFormula{b1} \left[ \tn{S} \vc{n} \right]_{\rm tan} +
\zeta \left[ \vu - \vc{V} \right]_{\rm tan}|_{\Gamma_\tau} = 0, \
\zeta \geq 0,
\eF
where $\zeta$ represents a ``friction'' coefficient.
If
$\zeta= 0$, we obtain the \emph{complete slip} while the
asymptotic limit $\zeta \to \infty$ gives rise to the standard
no-slip boundary conditions.

Concerning the heat flux we consider the conservative boundary
conditions
\bFormula{b2} \vq\cdot \vc{n}=0 \mbox { for all }
t\in [0,T],\  x \in  {\Gamma_t}.
\eF

For physical motivation of correct description of the fluid
boundary behavior, see Bul\' \i \v cek, M\' alek and Rajagopal
\cite{BMR1}, Priezjev and Troian \cite{PRTR} and the references
therein.

Finally, the problem \eqref{i2}-\eqref{b2} is supplemented by the
initial conditions
\bFormula{i7} \vr(0, \cdot) = \vr_0 ,\quad (\vr
\vu) (0, \cdot) = (\vr \vu)_0,\quad \vt(0,\cdot) = \vt_0,\quad  (\vr Q(\vt)) (0, \cdot) = (\vr Q)_0
= \vr_0 Q(\vt_0)
\quad \mbox{in}\ \Omega_0. \eF

The main goal of this paper is to show the \emph{existence} of global-in-time
\emph{weak solutions} to problem \eqref{i2}-\eqref{i7} for any
finite energy initial data. The existence theory for the
barotropic Navier-Stokes system on \emph{fixed} spatial domains in
the framework of weak solutions was developed in the seminal work
by Lions \cite{LI4}, and later extended in \cite{FNP} to a class
of physically relevant pressure-density state equations. Then it was extended by Feireisl to the full Navier-Stokes-Fourier system
\cite {EF70, EF71}.

The investigation of \emph{incompressible}
fluids in {time dependent} domains started with a seminal paper of
Ladyzhenskaya \cite{LAD2}, see also \cite{Neust1,NeuPen1,NeuPen2} for more recent results in this direction.

\emph{Compressible} fluid flows in time dependent domains in barotropic case
were examined in \cite{FeNeSt} for the \emph{no-slip} boundary conditions and in \cite{FKNNS} for the slip boundary conditions. The aim of this paper is to extend this result to the full system. 
We proceed in the following way.
\begin{enumerate}

\item In order to deal with the slip boundary condition, we introduce to the weak formulation of the momentum equation a term
\bFormula{i8}
\frac{1}{\ep} \int_0^T \int_{\Gamma_t}  (\vu - \vc{V} ) \cdot \vc{n} \ \vph \cdot \vc{n} \ {\rm dS}_{x} \ \dt,\ \ep > 0 \ \mbox{small},
\eF
which was originally proposed by Stokes and Carey in \cite{StoCar}. This extra term allows to consider the system in time independent domain which is divided by impermeable boundary $\Gamma_t$. In order to handle the behaviour of fluid in the solid domain $Q_T^c$ we use the following three (in fact four) level penalization scheme

\item In addition to (\ref{i8}), we introduce a \emph{variable}
shear viscosity coefficient $\mu = \mu_\omega$, where $\mu_\omega$
remains strictly positive in the fluid domain $Q_T$ but vanishes
in the solid domain $Q_T^c$ as $\omega \to 0$.
\item  We introduce the heat conductivity coefficient $\kappa_{\nu}(t,x,\vt)$ which
remains strictly positive in the fluid domain $Q_T$ but vanishes
in the solid domain $Q_T^c$ as $\nu \to 0$. For technical reasons this procedure needs to be done in two steps, more precisely we consider a step function $\kappa_\nu = \kappa$ in the fluid part and $\kappa_\nu = \nu\kappa$ in the solid part. This function is then mollified to get smooth $\kappa_{\nu,\xi}$, $\xi > 0$ with $\kappa_{\nu,\xi} \sil \kappa_\nu$ as $\xi \sil 0$.
\item Similarly to
the existence theory developed in \cite{FNP}, we introduce the
\emph{artificial pressure}
\[
p_\delta(\varrho,\vt ) = p(\varrho,\vt) + \delta \varrho^\beta,\
\beta \geq 4,\ \delta > 0,
\]
in the momentum equation (\ref{i3}).
\item Keeping $\ep, \delta, \nu, \xi$ and $\omega > 0$ fixed, we solve the modified problem in a (bounded) reference domain $B \subset \mathbb{R}^3$ chosen in such a way that
\[
\Ov{ \Omega}_\tau \subset B \ \mbox{for any}\ \tau \in [0,T].
\]
To this end, we adapt the existence theory for the compressible
Navier-Stokes system with variable viscosity coefficients
developed in \cite{EF71}.

\item We take the initial density
$\varrho_0$ vanishing outside $\Omega_0$ and letting $\ep \to 0$
for fixed $\delta , \ \omega > 0$ we obtain a ``two-fluid''
system where the density vanishes in the solid part $\left( (0,T)
\times B \right) \setminus Q_T$ of the reference domain together
with the thermal pressure $p_{\vt}$. It follows that the choice of the pressure is
fundamental.

\item Passing with $\xi \sil 0$ we recover the system with jump in the heat conductivity coefficient and justify the choice of the boundary condition of the test function in the weak formulation of the thermal energy balance.

\item Letting the viscosity vanish in the solid part, we perform the limit $\omega \to 0$, where the extra stresses disappear in the limit system. The desired conclusion results from the final limit process, where we set $\nu = \nu(\delta)$ and let $\delta \to 0$.

\end{enumerate}

The paper is organized as follows. In Section \ref{m}, we introduce all necessary preliminary material including a weak formulation of the problem and state the main result. Section \ref{p} is devoted to the penalized problem and to uniform bounds and existence of solutions at the starting level of approximations.
In Section \ref{s}, the limits for $\ep \to 0$, $\xi \sil 0$, $\omega \to 0$, and $\delta \to 0$ are preformed successively. Section \ref{d} discusses possible extensions and applications of the method.

\section{Preliminaries}
\label{m}

\subsection{Hypotheses}
\label{s:m1}

Hypotheses imposed on constitutive relations and
transport coefficients are motivated by the general existence
theory for the Navier-Stokes-Fourier system developed in
\cite{EF70,EF71}.

{\it Hypotheses on the pressure:}
%
%
\bFormula{hp1}\left.
\begin{array}{l}
p_e \in C[0,\infty)\cap C^1(0,\infty), \ p_e(0)=0,\\\\
p_e'(\vr) \geq a_1\vr^{\gamma -1}-b \mbox { for all } \vr >0,\\\\
p_e(\vr)\leq a_2\vr^{\gamma}+b \mbox { for all }\vr \geq 0,
\end{array}\right\}\eF
for certain constants $\gamma >\frac 32,$ 
$a_1, a_2, b>0$,
\bFormula{hp2} \left.\begin{array}{l} p_{\vt} \in C[0,\infty)\cap
C^1(0,\infty),\  p_{\vt} (0)=0,\\\\
p_{\vt} \mbox { is a non-decreasing function of } \vr \in
[0,\infty),\\\\
p_{\vt}(\vr)\leq c\vr ^{\frac{\gamma}{3} } \mbox { for all } \vr \geq
0.
\end{array}
\right\} \eF
\smallskip

{\it Hypotheses on the heat conductivity coefficient:}
\bFormula{hh1} \left.
\begin{array}{l}
\kappa = \kappa(\vt) \mbox { belongs to the class }
C^2[0,\infty),\\\\
k_1(\vt ^{\alpha }+1) \leq \kappa (\vt ) \leq k_2(\vt ^{\alpha }+1)
\mbox{ for all } \vt \geq 0,\\\\
\mbox{with constants }
k_1,\,k_2>0,\   \alpha \geq 4 \mbox { and } \alpha \geq \frac{12(\gamma - 1)}{\gamma}.
\end{array}
\right\} \eF
\smallskip

{\it Hypothesis on the thermal energy:}
\bFormula{hte1} \left.
\begin{array}{l}
Q=Q(\vt)= \int_0^{\vt}c_v(z) {\rm\, d}z,\\\\
\mbox{where }
c_v\in C^1[0,\infty) \mbox{ is such that there exist } 0<\underline c<\overline c<\infty , \\\\
\underline c(1+\vt^{\frac{\alpha}{2}-1})\leq c_v(\vt)\leq \overline c(1+\vt^{\frac{\alpha}{2}-1}).
\end{array}
\right\} \eF
We would like to emphasise that the lower bound on $c_v$ is restrictive (compare with \cite[Section 3]{EF70}). However, this assumption is essential since it allows us to derive estimates on the velocity  -- see \eqref{2.odhad.na.gradient}.


\subsection{Weak formulation, main result}

In the weak formulation, it is convenient to consider the continuity equation \eqref{i2} in the whole physical space $\mathbb{R}^3$ provided the density
$\vr$ is extended to be equal to zero outside the fluid domain, specifically
\bFormula{m1}
\int_{\Omega_\tau} \vr \varphi (\tau, \cdot) \ \dx - \int_{\Omega_0} \vr_0 \varphi (0, \cdot) \ \dx =
\int_0^\tau \int_{ \Omega_t} \left( \vr \partial_t \varphi + \vr \vu \cdot \Grad \varphi \right) \ \dxdt
\eF
for any $\tau \in [0,T]$ and any test function $\varphi \in \DC([0,T] \times \mathbb{R}^3)$. Moreover, equation (\ref{i2}) is also satisfied in the sense of
renormalized solutions introduced by DiPerna and Lions \cite{DL}:
\bFormula{m2}
\int_{\Omega_\tau} b(\vr) \varphi (\tau, \cdot) \ \dx - \int_{\Omega_0} b(\vr_0) \varphi (0, \cdot) \ \dx =
\int_0^\tau \int_{ \Omega_t} \left( b(\vr) \partial_t \varphi + b(\vr) \vu \cdot \Grad \varphi +
\left( b(\vr)  - b'(\vr) \vr \right) \Div \vu \varphi \right) \ \dxdt
\eF
for any $\tau \in [0,T]$, any $\varphi \in \DC([0,T] \times \mathbb{R}^3)$, and any $b \in C^1 [0, \infty)$, $b(0) = 0$, $b'(r) = 0$ for large $r$.
Of course, we suppose that $\vr \geq 0$ a.e. in $(0,T) \times \mathbb{R}^3$.

Similarly, the momentum equation \eqref{i3} is replaced by a family of integral identities
\bFormula{m3}
\int_{\Omega_\tau} \vr \vu \cdot \vph (\tau, \cdot) \ \dx - \int_{\Omega_0} (\vr \vu)_0 \cdot \vph (0, \cdot) \ \dx
\eF
\[
= \int_0^\tau \int_{\Omega_t} \left( \vr \vu \cdot \partial_t \vph + \vr [\vu \otimes \vu] : \Grad \vph + p(\vr,\vt) \Div \vph
- \tn{S} (\Grad \vu) : \Grad \vph\right) \dxdt
\]
for any $\tau \in [0,T]$ and any test function $\vph \in \DC([0,T] \times \mathbb{R}^3 ; \mathbb{R}^3)$ satisfying
\bFormula{m4}
\vph \cdot \vc{n}|_{\Gamma_\tau} = 0 \ \mbox{for any} \ \tau \in [0,T].
\eF

Then the impermeability condition (\ref{i6}) is satisfied in the sense of traces, specifically,
\bFormula{m5}
\vu,\nabla_x\vu \in L^2(Q_T; \mathbb{R}^3) \ \mbox{and}\ (\vu - \vc{V}) \cdot \vc{n}  (\tau , \cdot)|_{\Gamma_\tau}  = 0 \ \mbox{for a.a.}\ \tau \in [0,T].
\eF

According to \cite{EF70,EF71} the equation \eqref{i3a} is replaced by two inequalities, the \emph{thermal energy inequality}
\begin{equation}\label{i3a2}
\partial_t (\vr Q(\vt)) + \Div (\vr Q(\vt) \vu) - \Delta \K (\vt) \geq \tn{S}
: \Grad \vu - \vt p_{\vt}\Div \vu
\end{equation}
where
	$$
	\K(\vt) = \int_0^\vt \kappa (z) {\rm\,d}z
	$$
and the \emph{global total energy inequality}. In the case of the problem on a fixed domain $\Omega$, this inequality states simply
	\begin{equation}\label{i3a31}
	\int_{\Omega} \vr( \frac{1}{2} |\vu|^2 + P_e(\vr) + Q(\vt) )(\tau,\cdot) \,\dx
	\leq
	\int_{\Omega} \left( \frac{|(\vr\vu)_0|^2}{2\vr_0} + \vr_0 P_e(\vr_0) + \vr_0 Q(\vt_0) \right)(\cdot) \,\dx
	\end{equation}
for all $\tau \geq 0$, the elastic potential $P_e$ is defined in \eqref{p1pb}.

However, in the problem on the moving domain this inequality is no longer such simple and additional terms appear. We have
 \begin{equation}\label{i3a3}
	\int_{\Omega_\tau} \vr( \frac{1}{2} |\vu|^2 + P_e(\vr) + Q(\vt) )(\tau,\cdot) \,\dx
	\leq
	\int_{\Omega_0} \left( \frac{|(\vr\vu)_0|^2}{2\vr_0} + \vr_0 P_e(\vr_0) + \vr_0 Q(\vt_0)\right)(\cdot) \,\dx
	\end{equation}
	 \[
+ \int_{\Omega_\tau} (\vr \vu \cdot \vc{V}) (\tau, \cdot) \,\dx - \int_{\Omega_0} (\vr \vu)_{0} \cdot \vc{V}(0, \cdot) \,\dx
\]
\[
+ \int_0^\tau \int_{\Omega_t} \left(\tn{S}(\Grad \vu) : \Grad \vc{V} - \vr \vu \cdot \partial_t \vc{V} - \vr \vu \otimes \vu : \Grad \vc{V} -
p(\vr) \Div \vc{V} \right) \, \dxdt
\]
for all $\tau \geq 0$. We also emphasise that the boundary condition \eqref{b2} is replaced by the inequality
\begin{equation}
\frac{\partial \vt}{\partial \vc{n}} \geq 0 \mbox{ on }\Gamma_t \mbox{ for all } t\in [0,T].
\end{equation}

Following these considerations, in the weak formulation of the problem the temperature $\vt$ satisfies
	\bFormula{m6}
	\int_0^T\int_{\Omega_t} \vr Q(\vt) \partial_t \varphi + \vr Q(\vt) \vu \cdot \nabla \varphi +
	 \K(\vt) \Delta \varphi \,\dxdt
	\leq
	\int_0^T\int_{\Omega_t} (\vt p_\vt \Div \vu - \tn{S} : \Grad \vu ) \varphi \,\dxdt + \int_{\Omega_0} \vr_0Q(\vt_0)\varphi(0)
	\,\dx
	\eF
for any $\varphi \in C^\infty([0,T] \times \mathbb{R}^3),$ $\varphi \geq 0,$ $\varphi(T)=0$, $\nabla \varphi \cdot \vc{n}|_{\Gamma_\tau}=0$ for any $\tau \in [0,T]$.

\begin{Definition}\label{d:WS}
We say that the trio $(\vr,\vu,\vt)$ is a \emph{variational solution} of problem \eqref{i2}-\eqref{i3a} with boundary conditions \eqref{i6}-\eqref{b2} and initial conditions \eqref{i7} if
\begin{itemize}
\item $\vr \in L^\infty(0,T;L^\gamma(\mathbb{R}^3))$,
\item $\vu,\nabla_x\vu \in L^2(Q_T;\mathbb{R}^3)$, $\vr\vu \in L^\infty(0,T; L^{m}(\mathbb{R}^3; \mathbb{R}^3))$ for some $m > \frac 65$,
\item $\vr Q(\vt) \in L^\infty(0,T;L^1(\mathbb{R}^3)) \cap L^2(0,T;L^{q}(\mathbb{R}^3))$ for some $q > \frac 65$, $\log \vt, \vt p_\vt(\vr) \in L^2(Q_T)$, $\K(\vt) \in L^1(Q_T)$,
\item relations \eqref{m1}-\eqref{m5}, \eqref{i3a3} and \eqref{m6} are satisfied.
\end{itemize}
\end{Definition}


At this stage, we are ready to state the main result of the  present paper:

\bTheorem{m1}
Let $\Omega_0 \subset \mathbb{R}^3$ be a bounded domain of class $C^{2 + \nu}$, and let $\vc{V} \in C^1([0,T]; C^{3}_c (\mathbb{R}^3;\mathbb{R}^3))$ be given. Assume that the pressure $p(\vr,\vt)$ takes the form \eqref{p1p} with $p_e,p_\vt \in C[0, \infty) \cap C^1(0, \infty)$ complying with the hypothesis \eqref{hp1}, \eqref{hp2} with a certain $\gamma > 3/2.$ Moreover let $\vq$ be given by \eqref{q1} with $\kappa(\vt)$ complying with the hypothesis $\eqref{hh1}$ and let hypothesis  \eqref{hte1} 
be satisfied with $\alpha \geq \frac{12(\gamma - 1)} \gamma$.
%
Furthermore let the initial data fulfill
\[
\vr_0 \in L^\gamma (\mathbb{R}^3),\ \vr_0 \geq 0, \ \vr_0 \not\equiv 0,\ \vr_0|_{\mathbb{R}^3 \setminus \Omega_0} = 0,\
(\vr \vu)_0 = 0 \ \mbox{a.a. on the set} \ \{ \vr_0 = 0 \} ,\ \int_{\Omega_0} \frac{1}{\vr_0} |(\vr \vu)_0 |^2 \ \dx < \infty
\]
and
\[
\vt_0  \in L^\infty(\Omega_0),\quad \vt_0 \geq \underline{\vt}  > 0 \quad \mbox{ on } \Omega_0.
\]
Then the problem \eqref{i2}-\eqref{i3a} with boundary conditions \eqref{i6}-\eqref{b2} and initial conditions \eqref{i7} admits a variational solution on any time interval $(0,T)$ in the sense specified through Definition \ref{d:WS}.
	
	
\eT

The rest of the paper is devoted to the proof of Theorem \ref{Tm1}.

\section{Penalization}
\label{p}

For the sake of simplicity,
we restrict ourselves to the case $\zeta = 0$ in \eqref{b1} and $\eta = 0$ in \eqref{i4}. As we shall see in Section \ref{d}, the main ideas of the proof presented below require only straightforward modifications to accommodate the general case.

\subsection{Penalized problem - weak formulation}

Choosing $R > 0$ such that
$$
\vc{V} |_{[0,T] \times \{ |x| > R \} } = 0 ,\ \ \Ov{ \Omega }_0 \subset \{ |x| < R \}
$$
we take the reference domain $B = \{ |x| < 2R \}$.

Next, the shear viscosity coefficient  $\mu_\omega$ is taken such that
\bFormula{p2}
\mu_\omega \in \DC ([0,T] \times \mathbb{R}^3),\ 0 < \underline{\mu}_\omega \leq \mu_\omega (t,x) \leq \mu \ \mbox{in}\ [0,T] \times B, \
\mu_{\omega}(\tau, \cdot)|_{ \Omega_\tau } = \mu \ \mbox{for any} \ \ \tau \in [0,T]
\eF
and
\bFormula{p2_2}
 \mu_\omega \to 0 \mbox{ a.e. in } (0,T) \times B \setminus Q_T \qquad \mbox{ as } \omega \sil 0.
\eF

Similarly, we introduce variable heat conductivity coefficient as follows: $\kappa_\nu(t,x,\vt) = \chi_{\nu}(t,x) \kappa(\vt)$, where
\[
\chi_\nu = 1 \quad \mbox{ in } Q_T \qquad \mbox{ and } \chi_\nu = \nu \quad \mbox{ in } (0,T) \times B \setminus Q_T ,
\]
and its mollification $\kappa_{\nu,\xi}(t,x,\vt) = \chi_{\nu,\xi}(t,x)\kappa(\vt)$, where
\begin{equation}\label{chi}
\chi_{\nu,\xi} = \chi_\nu \star \Phi_\xi.
\end{equation}
Here $\star$ denotes the convolution in time-space (i.e. $\mathbb{R}^{4}$), $\Phi(t,x)$ is a standard mollifier such that $\Phi_\xi(t,x) \sil \delta_{t,x}$ as $\xi \sil 0$ in the sense of distributions.


Finally, let $\vr_0$, $(\vr\vu)_0$ and $\vt_0$ be initial conditions as specified in Theorem \ref{Tm1}. We define modified initial data $\vr_{0,\delta}$, $(\vr\vu)_{0,\delta}$ and $\vt_{0,\delta}$ so that 
\bFormula{data1}
\vr_{0, \delta} \geq 0, \ \vr_{0,\delta} \not\equiv 0,\ \vr_{0, \delta}|_{\mathbb{R}^3 \setminus \Omega_0} = 0,\ \int_{B}
\left( \vr_{0, \delta}^\gamma + \delta \vr_{0, \delta}^\beta \right) \dx \leq c,\ \vr_{0,\delta}\sil \vr_0 \ \mbox{in } L^\gamma(B), \ |\{\vr_{0,\delta}<\vr_0\}|\sil 0,
\eF
\bFormula{data2}
(\vr\vu)_{0,\delta} = \left\{
\begin{array}{ll}
(\vr\vu)_{0}\ & \mbox{if}\ \vr_{0,\delta}\geq \vr_0,\\
0 & \mbox{otherwise}
\end{array}
\right.
\eF
and for $\vt_{0,\delta}\in C^{2+\nu}(\mathbb R^3)$ it holds
	\bFormula{data3}
	\quad \nabla \vt_{0,\delta} \cdot \vc{n}|_{\Gamma_0} = 0, \  0< \underline{\vt} \leq \vt_{0,\delta} \leq \overline{\vt} \mbox{ on }\mathbb{R}^3 \mbox{ and } \vt_{0,\delta}|_{\Omega_0}\rightarrow \vt_0 \mbox{ in }L^1(\Omega_0).
	\eF

Now we are ready to state the weak formulation of the \emph{penalized problem}. Let $\beta > \max\{4,\gamma\}$.

 Again, we consider $\vr,\vu$ to be prolonged by zero outside of $(0,T) \times B$. The weak formulation of the continuity equations reads as
\bFormula{p3}
\int_{B} \vr \varphi (\tau, \cdot) \, \dx - \int_{B} \vr_{0,\delta} \varphi (0, \cdot) \, \dx =
\int_0^\tau \int_{B} \left( \vr \partial_t \varphi + \vr \vu \cdot \Grad \varphi \right) \,\dxdt
\eF
for any $\tau \in [0,T]$ and any test function $\varphi \in \DC([0,T] \times \mathbb{R}^3)$. The momentum equation is represented by the family of integral identities
\bFormula{p4}
\int_{B} \vr \vu \cdot \vph (\tau, \cdot) \, \dx - \int_{B} (\vr \vu)_{0,\delta} \cdot \vph (0, \cdot) \, \dx - \frac{1}{\ep} \int_0^\tau \int_{ \Gamma_t } \left( (\vc{V} - \vu ) \cdot \vc{n} \ \vph \cdot \vc{n} \right) \, {\rm dS}_{x} \, \dt
\eF
\[
= \int_0^\tau \int_{B} \left( \vr \vu \cdot \partial_t \vph + \vr [\vu \otimes \vu] : \Grad \vph + p(\vr,\vt) \Div \vph + \delta \vr^\beta \Div \vph
- \tn{S}_\omega : \Grad \vph \right) \, \dxdt
,
\]
\[
\mbox{ with }\quad \tn{S}_\omega = {\mu_\omega} \left( \Grad \vu + \Grad^t \vu - \frac{2}{3} \Div \vu \tn{I} \right)
\]
for any $\tau \in [0,T]$ and any test function $\vph \in \DC([0,T] \times B ; \mathbb{R}^3)$.

The thermal energy inequality in the penalized problem is satisfied in the renormalized sense:
	\bFormula{p42}
	\int_0^T\int_{B} (\vr+\delta) Q_h(\vt) \partial_t \varphi + \vr Q_h(\vt) \vu \cdot \nabla \varphi +
	 \K_{\nu,\xi,h}(t,x,\vt) \Delta \varphi +  \K_{h}(\vt) \nabla \chi_{\nu,\xi} \cdot \nabla_x \varphi  - \delta h(\vt)\vt^{\alpha + 1}\varphi \ \dxdt
	 \eF
	 \[	\leq \int_0^T\int_{B} (h(\vt)\vt p_\vt \Div \vu - h(\vt)(1-\delta)\tn{S}_\omega : \Grad \vu + h'(\vt)\kappa_{\nu,\xi}(t,x,\vt)|\nabla_x\vt|^2 ) \varphi \ \dxdt + \int_{B} (\vr_{0,\delta} + \delta)Q_h(\vt_{0,\delta}) \varphi(0) \ \dx
	\]
for any $\varphi \in C^\infty([0,T] \times \mathbb{R}^3),$ $\varphi \geq 0,$ $\varphi(T)=0$, $\nabla \varphi \cdot \vc{n}|_{\partial B}=0$ for any $\tau \in [0,T]$ and any $h \in C^\infty([0,\infty))$ such that
\begin{equation}\label{eq:h}
h(0) > 0, \quad h \text{ nonincreasing on } [0,\infty), \quad \lim_{z\rightarrow \infty} h(z) = 0, \quad h''(z)h(z) \geq 2h'(z)^2 \quad \text{ for all } z \geq 0,
\end{equation}
where $Q_h$, $\K_{\nu,\xi,h}$ and $\K_{h}$  are defined by
	$$Q_h=\int_0^\vt Q'(z) h(z) {\rm\,d}z, \quad\quad \K_{\nu,\xi,h}=\chi_{\nu,\xi}\int_0^\vt \kappa (z) h(z) {\rm\,d}z, \quad\quad \K_{h}=\int_0^\vt \kappa (z) h(z) {\rm\,d}z.
	$$
Finally, the global total energy inequality reads as
	\bFormula{p43}
	\int_0^T\int_{B} (-\partial_t\psi)\left( \vr\frac{1}{2} |\vu|^2 + \vr P_e(\vr) + \frac{\delta}{\beta-1}\vr^\beta + (\vr + \delta)Q(\vt) \right) \,\dx\dt + \delta\int_0^T \psi\int_B \tn{S}_\omega:\Grad \vu + \vt^{\alpha+1} \,\dx \dt
	\eF
	\[ \leq \int_{B}  \frac{1}{2} \frac{|(\vr\vu)_{0,\delta}|^2}{\vr_{0,\delta}} + \vr_{0,\delta} P_e(\vr_{0,\delta}) + \frac{\delta}{\beta-1}\vr_{0,\delta}^\beta + (\vr_{0,\delta} + \delta)Q(\vt_{0,\delta}) \, \dx
	 + \frac{1}{\ep} \int_0^T \int_{ \Gamma_t } \psi\left( (\vc{V} - \vu ) \cdot \vc{n} \ \vu \cdot \vc{n} \right) \, {\rm dS}_{x} \, \dt
	\]
for every $\psi \in C^{\infty}([0,T])$ satisfying
$$
\psi(0) = 1, \ \ \psi(T) = 0, \ \ \partial_t\psi\leq 0.
$$

\begin{Definition}\label{d:WSpen}
Let $\ep, \delta, \nu, \xi$ and $\omega$ be positive parameters and let $\beta > \max\{4,\gamma\}$.
We say that a trio $(\vr,\vu,\vt)$ is a renormalized solution to the penalized problem with initial data \eqref{data1}-\eqref{data3} if
\begin{itemize}
\item $\vr \in L^\infty(0,T;L^\gamma(\mathbb{R}^3)) \cap L^\infty(0,T;L^\beta(\mathbb{R}^3))$,
\item $\vu \in L^2(0,T; W^{1,2}_0(B; \mathbb{R}^3))$, $\vr\vu \in L^\infty(0,T; L^{m}(\mathbb{R}^3; \mathbb{R}^3))$ for some $m > \frac 65$,
\item $\vr Q(\vt) \in L^\infty(0,T;L^1(\mathbb{R}^3)) \cap L^2(0,T;L^{q}(\mathbb{R}^3))$ for some $q > \frac 65$, $\log \vt, \vt p_\vt(\vr) \in L^2((0,T)\times B)$, $\K_h(\vt) \in L^1((0,T)\times B)$ for all $h$ as in \eqref{eq:h},
\item relations \eqref{p3}-\eqref{p42} and \eqref{p43} are satisfied.
\end{itemize}
\end{Definition}

The choice of the no-slip boundary condition $\vu|_{\partial B} = 0$ is not essential here. We have the following existence theorem concerning weak solutions to the penalized problem.

\bTheorem{m2}
Let $\vc{V} \in C^1([0,T]; C^{3}_c (\mathbb{R}^3;\mathbb{R}^3))$ be given. Assume that the pressure satisfies the constitutive equation \eqref{p1p} with $p_e, p_\vt \in C[0, \infty) \cap C^1(0, \infty)$  complying with hypothesis  \eqref{hp1}, \eqref{hp2} with certain $\gamma > 3/2$. Moreover let $\vq$ be given by \eqref{q1} with assumption $\eqref{hh1}$, let hypothesis \eqref{hte1} be satisfied and let the initial data fulfill \eqref{data1}, \eqref{data2} and \eqref{data3}.
%
Finally, let $\beta > \max\{4,\gamma\}$, and $\ep, \omega, \nu,\xi, \delta > 0$.

Then the penalized problem admits a renormalized solution on any time interval $(0,T)$ in the sense specified by Definition \ref{d:WSpen}.
%
%
\eT


The existence of global-in-time solutions to the penalized problem can be shown by means of the method developed in \cite{EF71} to handle the nonconstant viscosity coefficients, more precisely see Proposition 4.3 in \cite{EF71}. Indeed, for $\ep > 0$ fixed, the extra penalty term in (\ref{p4}) can be treated as a ``compact'' perturbation. The presence of terms involving $\chi_{\nu,\xi}(x,t)$ in the thermal energy inequality \eqref{p42} does not cause any additional difficulties since this function is smooth and bounded away from zero.


In addition, since $\beta > 4$, the density is square integrable and we may use the regularization technique of DiPerna and Lions \cite{DL}
to deduce the renormalized version of (\ref{p3}), namely
\bFormula{p6a}
\int_{B} b(\vr) \varphi (\tau, \cdot) \, \dx - \int_{B} b(\vr_{0,\delta}) \varphi (0, \cdot) \, \dx =
\int_0^\tau \int_{B} \left( b(\vr) \partial_t \varphi + b(\vr) \vu \cdot \Grad \varphi +
\left( b(\vr)  - b'(\vr) \vr \right) \Div \vu \varphi \right) \,\dxdt
\eF
for any $\varphi$ and $b$ as in (\ref{m2}).


\subsection{Modified energy inequality and uniform bounds}\label{s:bounds}

Since the vector field $\vc{V}$ vanishes on the boundary $\partial B$ it may be used as a test function in (\ref{p4}).  Combining the resulting expression with the energy inequality
(\ref{p43}), we obtain
\bFormula{p7}
\int_B \left( \frac{1}{2} \vr |\vu|^2 + \vr P_e (\vr) + \frac{\delta}{\beta - 1} \vr^\beta + (\vr + \delta) Q(\vt) \right)(\tau, \cdot) \, \dx
\eF
\[
+ \  \delta\int_0^\tau \int_B \tn{S}_\omega:\Grad \vu + \vt^{\alpha+1} \ \dx \dt + \frac{1}{\ep} \int_0^\tau \int_{\Gamma_t} \left| \left( \vu - \vc{V} \right) \cdot \vc{n} \right|^2 \ {\rm dS}_x \, \dt
\]
\[
\leq \int_B \left( \frac{1}{2 \vr_{0,\delta}} |(\vr \vu)_{0,\delta} |^2 + \vr_{0,\delta} P_e(\vr_{0, \delta}) + \frac{\delta}{\beta - 1} \vr_{0, \delta}^\beta  + (\vr_{0,\delta} + \delta) Q(\vt_{0,\delta})\right) \, \dx
\]
\[
+ \int_B \Big(  (\vr \vu \cdot \vc{V}) (\tau, \cdot) - (\vr \vu)_{0,\delta} \cdot \vc{V}(0, \cdot) \Big)
\, \dx
\]
\[
+ \int_0^\tau \int_B \left( \mu_\omega \left(\Grad \vu + \Grad^t \vu - \frac{2}{3} \Div \vu \tn{I} \right) : \Grad \vc{V} - \vr \vu \cdot \partial_t \vc{V} - \vr \vu \otimes \vu : \Grad \vc{V} -
p(\vr,\vt) \Div \vc{V} - \frac{ \delta }{\beta - 1} \vr^\beta \Div \vc{V} \right)  \dxdt.
\]
Since the vector field $\vc{V}$ is regular suitable manipulations with the H\"older, Young and Poincare inequalities and thermodynamical hypothesis yield
\bFormula{p7a}
\int_B \left( \frac{1}{2} \vr |\vu|^2 + \vr P_e (\vr) + \frac{\delta}{\beta - 1} \vr^\beta + (\vr + \delta) Q(\vt) \right)(\tau, \cdot) \, \dx
+ \frac{\delta}{2}\int_0^\tau \int_B \tn{S}_\omega:\Grad \vu  \ +\  \vt^{\alpha+1} \,\dx \dt
\eF
\[
+ \frac{1}{\ep} \int_0^\tau \int_{\Gamma_t} \left| \left( \vu - \vc{V} \right) \cdot \vc{n} \right|^2 \ {\rm dS}_x \, \dt
\]
\[
\leq \int_B \left( \frac{1}{2 \vr_{0,\delta}} |(\vr \vu)_{0,\delta} |^2 + \vr_{0,\delta} P_e(\vr_{0, \delta}) + \frac{\delta}{\beta - 1} \vr_{0, \delta}^\beta  + (\vr_{0,\delta} + \delta) Q(\vt_{0,\delta})\right) \, \dx  - \int_B  (\vr \vu)_{0,\delta} \cdot \vc{V}(0, \cdot) \,\dx + C(\vc{V}, \delta ).
\]

As
\[
p_e(\vr) \leq c(1 + \vr P_e(\vr)) \ \mbox{for all} \ \vr \geq 0,
\]
relation (\ref{p7a}) gives rise to the following bounds
independent of parameters $\ep$ and $\xi$:
\bFormula{p8}
{\rm ess} \sup_{t \in (0,T)} \| \sqrt{\vr} \vu (t , \cdot)  \|_{L^2 (B;\mathbb{R}^3)} \leq c,
\eF
\bFormula{p9}
{\rm ess} \sup_{t \in (0,T)} \int_B \vr P_e(\vr)(t , \cdot) \, \dx \leq c \quad \mbox{yielding} \quad {\rm ess} \sup_{t \in (0,T)} \left\| \vr (t, \cdot) \right\|_{L^\gamma(B)} \leq c,
\eF
\bFormula{p10}
{\rm ess} \sup_{t \in (0,T)} \delta \| \vr (t, \cdot) \|^\beta_{L^\beta(B)} \leq c.
\eF
Since
	\begin{equation*}
	\frac{ \mu_{\omega} }{2} \left| \Grad \vu + \Grad^t \vu - \frac{2}{3} \Div \vu \tn{I} \right|^2 =
	\left[  \mu_\omega \left( \Grad \vu + \Grad^t \vu \right) - \frac{2 \mu_\omega }{3} \Div \vu \tn{I} \right] : \Grad \vu ,
	\end{equation*}
we obtain that
\bFormula{p11}
\int_0^T \int_B \mu_{\omega} \left| \Grad \vu + \Grad^t \vu - \frac{2}{3} \Div \vu \tn{I} \right|^2 \, \dxdt \leq c.
\eF
Moreover,
\bFormula{p12}
\int_0^T \int_{\Gamma_t} \left| (\vu - \vc{V}) \cdot \vc{n}  \right|^2 \ {\rm dS}_{x} \ \dt  \leq \ep c.
\eF
We note also that the total mass is conserved, meaning
\bFormula{p13}
\int_B \vr (\tau, \cdot) \, \dx = \int_B \vr_{0, \delta} \, \dx = \int_{\Omega_0} \vr_{0, \delta} \, \dx \leq c \ \mbox{for any}\ \tau \in [0,T].
\eF
Thus, relations (\ref{p8}), (\ref{p11}), (\ref{p13}), combined with the generalized version of Korn's inequality (see \cite[Theorem 10.17]{FeNo6}),
imply that
\bFormula{p14}
\int_0^T \| \vu (t, \cdot) \|^2_{W^{1,2}_0 (B;\mathbb{R}^3)} \,\dt \leq c(\omega).
\eF
 Moreover directly from the energy estimates we get that
 	\bFormula{p142}
	{\rm ess}\sup_{t\in(0,T)} \| (\vr + \delta)  Q(\vt) \|_{L^1(B)} \leq c,
	\eF
	\bFormula{p145}
	\| \vt \|_{L^{\alpha + 1}((0,T) \times B)} \leq c .
	\eF
By the proper choice of renormalization function  $h$  the same methods as in Section 5.1 \cite{EF71}  gives
	\bFormula{p143}
	\| \delta  \nabla \log \vt \|_{L^2((0,T) \times B)} \leq c,
	\eF
	\bFormula{p144}
	\| \delta \nabla \vt^{\alpha/2} \|_{L^2((0,T) \times B)} \leq c,
	\eF
and moreover
	\bFormula{p146}
	\|\delta \nabla \vt^{\frac{\alpha + 1 - \lambda}{2}} \|_{L^2((0,T) \times B)} \leq c \quad \mbox{ for any } 0< \lambda <1 .
	\eF
We emphasise that all constants in \eqref{p8} -- \eqref{p146} depend on $\delta$.

\subsection{Pressure estimates}

We use the technique based on the Bogovskii operator in order to derive the following estimate
\bFormula{p15}
\int\int_{K} \Big( \left(p_e(\vr)+ \vt p_\vt(\vr) \right) \vr^\lambda + \delta \vr^{\beta + \lambda} \Big) \, \dxdt \leq c({K}) \ \mbox{for a certain} \ \lambda > 0
\eF
for any compact $K \subset [0,T] \times \Ov{B}$ such that
\[
{K} \cap \left( \cup_{ \tau \in [0,T] } \Big( \{ \tau \} \times \Gamma_\tau \Big) \right) = \emptyset.
\]
For details we refer reader to \cite[Section 4.2]{EF71} or \cite{FP13}.



\section{Singular limits}

\label{s}

In this section, we perform successively the singular limits $\ep \to 0$, $\xi \sil 0$ $\omega \to 0$, and $\delta \to 0$ with $\nu=\nu(\delta)$.

\subsection{Penalization limit}\label{s:penlimit}

Firstly, we proceed with $\ep \to 0$ in \eqref{p3}, \eqref{p4}, \eqref{p42} and \eqref{p7} while other parameters $\xi$, $\nu$, $\omega$ and $\delta$ remain fixed.
Let $\{ \vre, \vue, \vte \}_{\ep > 0}$ be the corresponding sequence of renormalized solutions of the penalized problem given by Theorem \ref{Tm2}. The estimates \eqref{p9}, \eqref{p14} together with the equation of continuity \eqref{p3}, imply that
\[
\vre \to \vr \ \mbox{in}\ C_{\rm weak}([0,T] ; L^\gamma (B)),
\]
and, up to a subsequence,
\[
\vue \to \vu \ \mbox{weakly in}\ L^2(0,T; W^{1,2}_0 (B, \mathbb{R}^3)).
\]

Directly from \eqref{p12} we derive in the limit as $\ep \to 0$ that
\bFormula{s1}
( \vu - \vc{V} ) \cdot \vc{n} (\tau, \cdot) |_{\Gamma_\tau} = 0 \ \mbox{for a.a.}\ \tau \in [0,T].
\eF

Consequently, in accordance with (\ref{p8}), (\ref{p9}) and the compact embedding $L^\gamma(B) \hookrightarrow\hookrightarrow W^{-1,2}(B)$,
we obtain
\bFormula{s1a}
\vre \vue  \to \vr \vu \ \mbox{weakly-(*) in}\ L^\infty(0,T; L^{2 \gamma / (\gamma + 1)}(B;\mathbb{R}^3)),
\eF
and, due to the embedding $W^{1,2}_0 (B) \hookrightarrow L^6(B)$,
\[
\vre \vue \otimes \vue \to \Ov{ \vr \vu \otimes \vu } \ \mbox{weakly in}\ L^2(0,T; L^{6 \gamma / (4 \gamma + 3)} (B; \mathbb{R}^3)),
\]
where the bar denotes a weak limit of a composed function.

Finally we deduce from \eqref{p4} that
\[
\vre \vue \to \vr \vu \ \mbox{in} \ C_{\rm weak}([T_1, T_2]; L^{2 \gamma / (\gamma + 1)}(O;\mathbb{R}^3))
\]
for any space-time cylinder
\[
(T_1, T_2) \times {O} \subset [0,T] \times B,\  [T_1, T_2] \times \Ov{O} \cap \cup_{\tau \in [0,T]} \left( \{ \tau \} \times \Gamma_\tau \right) =
\emptyset.
\]

Since $L^{2 \gamma / (\gamma + 1)}(B) \hookrightarrow\hookrightarrow W^{-1,2}(B)$, we conclude that
\[
\Ov{\vr \vu \otimes \vu} = \vr \vu \otimes \vu \ \mbox{a.a. in}\ (0,T) \times B.
\]

\subsubsection{Strong convergence of temperature}

By assumption \eqref{hte1}
	\begin{multline}
	\|\nabla Q(\vte)\|_{L^2 ((0,T)\times B)}^2  \leq  C\left(\int_{(0,T)\times B} |\nabla \vte|^2 \, \dx \dt + \int_{(0,T)\times B}|\nabla \vte^{\frac \alpha 2}|^{2}\, \dx \dt\right)\\
	\leq C \left(\int_{\{\vte<1\}} \left|\frac{\nabla \vte}{\vte}\right|^2\, \dx \dt + \int_{\{\vte\geq 1\}} |\nabla \vte|^2|\vte|^{\alpha - 2}
	\, \dx \dt + \int_{(0,T)\times B}|\nabla \vte^{\frac \alpha 2}|^{2}\ \dx \dt\right)\\ \leq C(\delta)\left(\|\nabla \log \vte\|^2_{L^2((0,T)\times B)} + \|\nabla \vte^{\frac \alpha 2}\|^2_{L^2((0,T)\times B)} \right)
	\end{multline}
and, with help of \eqref{p143} and \eqref{p144}, we get
	\begin{equation}\label{s10b}
	\| Q(\vte)  \|_{L^2(0,T;W^{1,2}(B))} \leq c,
	\end{equation}
consequently,
	\begin{equation}\label{s10c}
	Q(\vte) \to \overline{Q(\vt)} \mbox{ weakly in }L^2(0,T;W^{1,2}(B)).
	\end{equation}
Next, according to Lemma~6.3 in \cite{EF70} applied to  $(\delta + \vre) Q(\vte)$ we obtain
	\begin{equation}\label{s10d}
	(\delta + \vre) Q(\vte) \to (\delta + \vr) \overline{Q(\vt)} \mbox{ strongly in }L^2(0,T; W^{-1,2} (B)).
	\end{equation}
Similarly as in Subsection 7.3.6 in \cite{EF70} one gets
\begin{equation}\label{s10e}
\vte \rightarrow \vt \mbox{ strongly in } L^2((0,T)\times \Omega).
\end{equation}

\subsubsection{Pointwise  convergence of density}\label{s:scd}

Next, we show pointwise convergence of the sequence $\{ \vre \}_{\ep > 0}$. 
We define
\[
p_\delta (\vr,\vt) = p(\vr,\vt) + \delta \vr^\beta ,\ T_k(\vr) = \min \{ \vr , k \}.
\]
Similarly as in \cite{EF71} we establish the effective viscouse pressure identity:
\bFormula{s5}
\Ov{ p_\delta (\vr,\vt) T_k (\vr) } - \Ov{ p_\delta (\vr,\vt) } \ \Ov{ T_k (\vr) } = \frac{4}{3} \mu_\omega \left(
\Ov{ T_k (\vr) \Div \vu } - \Ov{T_k (\vr) } \Div \vu \right),
\eF
which holds only on compact sets ${K} \subset [0,T] \times B$ satisfying
\[
{K} \cap \left( \cup_{ \tau \in [0,T] } \Big( \{ \tau \} \times \Gamma_\tau \Big) \right) = \emptyset.
\]

Following \cite{EF71}, we introduce the \emph{oscillations defect measure}
\[
{\bf osc}_q [ \vre \to \vr] ({K}) = \sup_{k \geq 0} \left( \limsup_{\ep \to 0}
\int_{{K}} | T_k (\vre) - T_k(\vr) |^q \, \dxdt \right),
\]
and use (\ref{s5}) to conclude that
\bFormula{s6}
{\bf osc}_{\gamma + 1} [ \vre \to \vr] ({K}) \leq c(\omega) < \infty,
\eF
where the constant $c$ is \emph{independent} of ${K}$. Thus
\bFormula{s7}
{\bf osc}_{\gamma + 1} [ \vre \to \vr] ([0,T] \times B) \leq c(\omega),
\eF
which implies, by virtue of the procedure developed in \cite{EF70}, the desired conclusion
\bFormula{s8}
\vre \to \vr \ \mbox{a.e. in}\ (0,T) \times B.
\eF

\subsubsection{Passing to the limit with $\varepsilon$}

Passing to the limit in \eqref{p3} we obtain
\bFormula{s2}
\int_{B} \vr \varphi (\tau, \cdot) \, \dx - \int_{B} \vr_{0,\delta} \varphi (0, \cdot) \, \dx =
\int_0^\tau \int_{B} \left( \vr \partial_t \varphi + \vr \vu \cdot \Grad \varphi \right) \, \dxdt
\eF
for any $\tau \in [0,T]$ and any $\varphi \in \DC([0,T] \times \mathbb{R}^3)$.

The limit in the momentum equation \eqref{p4} is more delicate. Since we have at hand only the \emph{local estimates} (\ref{p15}) on the pressure,
we have to restrict ourselves to the class of test functions
\bFormula{s4}
\vph \in C^1([0,T]; W^{1, \infty}_0 (B; \mathbb{R}^3)),\ {\rm supp}[ \Div \vph (\tau, \cdot)] \cap \Gamma_\tau = \emptyset,\
\vph \cdot \vc{n}|_{\Gamma_\tau} = 0 \ \mbox{for all}\ \tau \in [0,T].
\eF


In accordance with \eqref{s8} and \eqref{p15}, the momentum equation reads
\bFormula{s9}
\int_{B} \vr \vu \cdot \vph (\tau, \cdot) \, \dx - \int_{B} (\vr \vu)_{0,\delta} \cdot \vph (0, \cdot) \, \dx
\eF
\[
= \int_0^\tau \int_{B} \left( \vr \vu \cdot \partial_t \vph + \vr [\vu \otimes \vu] : \Grad \vph + {p(\vr,\vt)} \Div \vph + \delta {\vr^\beta} \Div \vph
- \mu_\omega \left( \Grad \vu + \Grad^t \vu - \frac{2}{3} \Div \vu \tn{I} \right) : \Grad \vph \right) \, \dxdt
\]
for any test function $\vph$ as in (\ref{s4}). In addition, as already observed, the limit solution $\{ \vr, \vu \}$ satisfies also the renormalized equation (\ref{p6a}).
The strong convergence of the temperature \eqref{s10e} is sufficient for the non-linear terms in thermal energy equation to pass to their limits counterparts.
In particular:
{	$$
	\vre Q(\vte) \to \vr Q (\vt) \mbox{ weakly in } L^2(0,T;L^{q_1}(B)) \mbox{ with }q_1=\frac{6\gamma}{6 + \gamma}
	$$
and
	$$
	\vre Q(\vte) \vue \to \vr Q (\vt) \vu \mbox{ weakly in } L^2(0,T;L^{q_2}(B)) \mbox{ with }q_2=\frac{6\gamma}{3+4\gamma} .
	$$
	}
Moreover, if $p_\vt$ satisfies hypothesis \eqref{hp2} we have
	$$
	\vte p_\vt(\vre) \Div \vue \to \vt p_\vt(\vr) \Div \vu \mbox{ weakly in } L^1((0,T) \times B)	
	$$
(in order to provide it one may use \eqref{p144} and additional information from artificial pressure)
and since $h$ satisfies \eqref{eq:h},
	$$
	h(\vte) \vte p_\vt(\vre) \Div \vue \to h(\vt) \vt p_\vt(\vr) \Div \vu \mbox{ weakly in } L^1((0,T) \times B).	
	$$	
Due to the convexity of the function
	\begin{equation*}
	[\tn{M},\vt] \mapsto
		\left\{
		\begin{array}{cc}
	      h(\vt) \left( \frac{\mu}{2} \tn{M} : \tn{M}  + \lambda (\mbox{tr} [ \tn{M} ])^2 \right) \mbox{ if } \vt \geq 0,\, \tn{M} \in \mathbb R^{3\times 3}  \\
	      \infty \mbox{ if } \vt <0 \\
	   \end{array}
	   \right.
	\end{equation*}
and \cite[Lemma~4.8]{EF70}, we get
	\begin{equation}
	\int_0^T \int_{B} h({\vt}) \tn{S}_\omega( \Grad \vu) :\Grad \vu \varphi \, \dxdt \leq \liminf\limits_{\ep \to 0}  \int_0^T \int_{B} h({\vte}) \tn{S}_\omega(\Grad \vue) :\Grad \vue \varphi \, \dxdt
	\end{equation}
where $\varphi$ is as in \eqref{p42}.
In order to treat the most-right term in \eqref{p42} we apply \cite[Corollary~2.2]{EF70}. Let us set
	\begin{equation}\label{s11}
	\int_0^T \int_{B} \chi_{\nu,\xi}\kappa(\vte)(- h'(\vte)) |\nabla \vte|^2 \,\dxdt = \int_0^T \int_{B} |\sqrt{\chi_{\nu,\xi}} G'(\vte)\nabla \vte|^2 \,\dxdt
	\quad \mbox{ with } G'(z) = \sqrt{ \kappa(z)(- h'(z))}.
	\end{equation}
The sequence $\chi_{\nu,\xi} \nabla G(\vte)$ is uniformly bounded in $L^{1+ \lambda}((0,T)\times B)$ for certain $\lambda >0$ and, consequently, $\sqrt{\chi_{\nu,\xi}} \nabla G(\vte) \to \sqrt{\chi_{\nu,\xi}} \nabla \overline{G(\vt)}$ weakly in $L^1((0,T)\times B)$. Due to a.e. convergence of $\vte \to \vt$ we have $\sqrt{\chi_{\nu,\xi}} \sqrt{\kappa(\vte)(-h'(\vte))}\to \sqrt{\chi_{\nu,\xi}} \sqrt{\kappa(\vt)(-h'(\vt))}$ strongly in $L^2$. Thus, as $\nabla\vte\to\nabla\vt$ weakly in $L^2$, we get $\sqrt{\chi_{\nu,\xi}} \nabla G(\vte)\to \sqrt{\chi_{\nu,\xi}} \nabla G(\vt)$ weakly in $L^1((0,T)\times B)$. Then convexity of $\Phi(\cdot ) = | \cdot |^2$ and \cite[Corollary~2.2]{EF70} provide
 	\begin{equation}\label{s12}
	\int_0^T \int_{B} \chi_{\nu,\xi}\kappa(\vt)(- h'(\vt)) |\nabla \vt|^2 \varphi \,\dxdt \leq \liminf\limits_{\ep \to 0} \int_0^T \int_{B}\chi_{\nu,\xi}\kappa(\vte)(- h'(\vte)) |\nabla \vte|^2 \varphi\, \dxdt ,
	\end{equation}
where $\varphi$ is as in \eqref{p42}. Since $\vte$ is uniformly bounded in $L^{\alpha + 1}$ and due to properties of $h$, the sequence $h(\vte)\vte^{\alpha + 1}$ is uniformly integrable in $L^1$ (see \cite[Proposition~2.1]{EF70}). Then a.e. convergence of $\vte$ provides
 	\begin{equation*}
	\int_0^T \int_B h(\vte) \vte^{\alpha+1} \varphi \, \dxdt \to \int_0^T \int_B h(\vt) \vt^{\alpha+1} \varphi \, \dxdt  \quad \mbox{ as } \ep \to 0.
	\end{equation*}
By the same token
	\begin{equation*}
	\int_0^T \int_B \K_{\nu,\xi,h} (\vte) \Delta \varphi \, \dxdt \to \int_0^T \int_B \K_{\nu,\xi,h}(\vt) \Delta \varphi \,\dxdt  \quad \mbox{ as } \ep \to 0,
	\end{equation*}
and also
\begin{equation*}
	\int_0^T \int_B \K_{h} (\vte) \nabla \chi_{\nu,\xi} \cdot \nabla \varphi \, \dxdt \to \int_0^T \int_B \K_{h}(\vt) \nabla\chi_{\nu,\xi} \cdot \nabla \varphi \,\dxdt  \quad \mbox{ as } \ep \to 0,
	\end{equation*}
where $\varphi$ is as in \eqref{p42}.
Thus, the renormalized thermal energy inequality \eqref{p42} takes in the limit the same form. 

Finally, we obtain the following energy inequality when passing $\ep \to 0$
	\begin{multline}\label{s43}
	\int_{B} \left( \vr\frac{1}{2} |\vu|^2 + \vr P_e(\vr) + \frac{\delta}{\beta-1}\vr^\beta + (\vr + \delta)Q(\vt) \right)(\tau) \,\dx + \delta\int_0^\tau \int_B \tn{S}_\omega:\Grad \vu + \vt^{\alpha+1} \,\dx \dt
\\ \leq \int_{B}  \frac{1}{2} \frac{|(\vr\vu)_{0,\delta}|^2}{\vr_{0,\delta}} + \vr_{0,\delta} P_e(\vr_{0,\delta}) + \frac{\delta}{\beta-1}\vr_{0,\delta}^\beta + (\vr_{0,\delta} + \delta)Q(\vt_{0,\delta}) \,\dx
+ \int_B \Big(  (\vr \vu \cdot \vc{V}) (\tau, \cdot) - (\vr \vu)_{0,\delta} \cdot \vc{V}(0, \cdot) \Big)
\ \dx\\
+ \int_0^\tau \int_B \left( \mu_\omega \left(\Grad \vu + \Grad^t \vu - \frac{2}{3} \Div \vu \tn{I} \right) : \Grad \vc{V} - \vr \vu \cdot \partial_t \vc{V} - \vr \vu \otimes \vu : \Grad \vc{V} -
p(\vr,\vt) \Div \vc{V} - \frac{ \delta }{\beta - 1} \vr^\beta \Div \vc{V} \right) \, \dxdt.
\end{multline}

\subsubsection{Fundamental lemma and extending the class of test functions}
\label{a}
Our next goal is to use the specific choice of the initial data $\vr_{0, \delta}$ to get rid of the density-dependent terms in (\ref{s9}) supported by the ``solid'' part $\left( (0,T) \times B \right) \setminus Q_T$. To this end, we proceed the same way as in the barotropic case, more precisely we use \cite[Lemma 4.1]{FKNNS}.

\bLemma{a1}
Let $\vr \in L^\infty (0,T; L^2(B))$, $\vr \geq 0$,  $\vu \in L^2(0,T; W^{1,2}_0(B;\mathbb{R}^3))$ be a weak solution of the equation of continuity,
specifically,
\bFormula{a1}
\int_{B} \Big( \vr (\tau, \cdot) \varphi (\tau, \cdot) - \vr_0 \varphi(0, \cdot) \Big) \,\dx
= \int_0^\tau \int_{B} \Big( \vr \partial_t \varphi + \vr \vu \cdot \Grad \varphi \Big) \, \dx\dt
\eF
for any $\tau \in [0,T]$ and any test function $\varphi \in C^1_c ([0,T] \times \mathbb{R}^3)$.

In addition, assume that
\bFormula{a2}
( \vu - \vc{V} )(\tau, \cdot) \cdot \vc{n} |_{\Gamma_\tau}  = 0 \ \mbox{for a.a.}\ \tau \in (0,T),
\eF
and that
\[
\vr_0 \in L^2 (\mathbb{R}^3), \ \vr_0 \geq 0,  \ \vr_0 |_{B \setminus \Omega_0} = 0.
\]

Then
\[
\vr(\tau, \cdot) |_{B \setminus \Omega_\tau} = 0 \ \mbox{for any}\ \tau \in [0,T].
\]

\eL

\bigskip

By virtue of Lemma \ref{La1}, the momentum equation (\ref{s9}) reduces to
\bFormula{s20}
\int_{\Omega_\tau } \vr \vu \cdot \vph (\tau, \cdot) \, \dx - \int_{\Omega_0} (\vr \vu)_{0, \delta} \cdot \vph (0, \cdot) \, \dx
\eF
\[
= \int_0^\tau \int_{\Omega_t} \left( \vr \vu \cdot \partial_t \vph + \vr [\vu \otimes \vu] : \Grad \vph + {p(\vr,\vt)} \Div \vph + \delta {\vr^\beta} \Div \vph
- \mu_\omega \left( \Grad \vu + \Grad^t \vu - \frac{2}{3} \Div \vu \tn{I} \right) : \Grad \vph \right) \, \dxdt
\]
\[
- \int_0^\tau \int_{ B \setminus \Omega_t } \mu_\omega \left( \Grad \vu + \Grad^t \vu - \frac{2}{3} \Div \vu \tn{I} \right) : \Grad \vph  \, \dxdt
\]
for any test function $\vph$ as in (\ref{s4}). We remark that in this step we crucially need the extra pressure term $\delta \vr^\beta$ ensuring the density $\vr$ to be square integrable.

Next, we argue the same way as in \cite[Section 4.3.1]{FKNNS} to conclude, that the momentum equation \eqref{s20} holds in fact for any test function $\vph$ such that
\begin{equation}\label{s4a}
\vph \in C^\infty_c([0,T]\times \mathbb{R}^3;\mathbb{R}^3), \qquad \vph(\tau, \cdot)\cdot \vc{n} |_{\Gamma_\tau} = 0 \quad \text{ for any } \tau \in [0,T].
\end{equation}

\subsection{The discontinuous heat conductivity coefficient}

The aim of this section is to proceed to a limit with $\xi \sil 0$. In other words we pass from a smooth function $\chi_{\nu,\xi}$ presented in \eqref{chi} to a jump function $\chi_\nu$. We denote by $\vu_\xi$, $\vt_\xi$ and $\vr_\xi$ the sequences of functions satisfying relations \eqref{s2}, \eqref{s20}, \eqref{p42} and \eqref{s43} which were constructed in a previous section. To pass to the limit $\xi \to 0$ we use the a priori estimates from Section \ref{s:bounds} to conclude that there exist $\vr$, $\vu$ and $\vt$ such that
\begin{equation}\label{conv.jump.k}
\begin{split}
\vr_\xi\to \vr &\mbox{ in } C_{\rm weak}([0,T];L^{\gamma} (B)) \cap C_{\rm weak}([0,T];L^{\beta} (B)),\\
\vu_\xi\to \vu & \mbox{ weakly in } L^2(0,T;W^{1,2}_0(B, \mathbb{R}^3)),\\
\vt_\xi \to \vt &\mbox{ strongly in }L^2((0,T)\times B).
\end{split}
\end{equation}
Using the same procedures as in Section \ref{s:penlimit} we pass to the limit in most of the terms in the equations and inequalities \eqref{s2}, \eqref{s20}, \eqref{p42} and \eqref{s43}. We treat in detail only the terms directly involving $\chi_{\nu,\xi}$.

Firstly, we consider the term
$$
I_\xi:=\int_0^T\int_B \chi_{\nu,\xi}\K_h(\vt_\xi) \Delta \varphi \, \dxdt + \K_h(\vt_\xi) \nabla \chi_{\nu,\xi} \cdot \nabla \varphi
\, \dxdt,$$
with $\varphi$ as for \eqref{p42}.
By integration by parts we derive
$$
-I_\xi= \int_0^T\int_B \chi_{\nu,\xi}\nabla \K_h(\vt_\xi) \cdot \nabla \varphi \ \dxdt= \int_0^T\int_B\chi_{\nu,\xi} \kappa(\vt_\xi)h(\vt_\xi) \nabla \vt_\xi \cdot \nabla \varphi \ \dxdt.
$$
As a corollary of \eqref{p142}, \eqref{p145}, \eqref{p146} and \eqref{conv.jump.k}$_3$ we get (up to a subsequence)
\begin{equation*}
\begin{split}
\vartheta_\xi &\rightarrow \vartheta \ \mbox{ weakly in }\ L^2(0,T; W^{1,2}(B)) ,\\
\vartheta_\xi &\rightarrow \vartheta \ \mbox{ strongly in }\ L^p(0,T; L^p(B)) \ \mbox{ for all }\ p \in [1, \alpha +1) , \\
\vartheta^{\frac \alpha 2 + \frac 5{12}}_\xi &\rightarrow \vartheta^{\frac \alpha 2+ \frac 5{12}} \ \mbox{ weakly in }\ L^2(0,T; W^{1,2}(B)) ,
\end{split}
\end{equation*}
and also
\begin{equation}\label{convergence1}
\chi_{\nu,\xi}(t,x) \frac{\kappa(\vartheta_\xi)h(\vartheta_\xi)}{1+\vartheta_\xi^{\frac \alpha 2 - \frac 7{12}}} \rightarrow \chi_{\nu}(t,x)\frac{\kappa(\vartheta)h(\vartheta)}{1+\vartheta^{\frac \alpha 2 - \frac{7}{12}}}\ \mbox{ strongly in } L^2((0,T)\times B).
\end{equation}
Indeed, if we denote $f_\xi = \chi_{\nu,\xi}(t,x) \frac{\kappa(\vartheta_\xi)h(\vartheta_\xi)}{1+\vartheta_\xi^{\frac \alpha 2 - \frac 7{12}}}$ and $f = \chi_{\nu}(t,x)\frac{\kappa(\vartheta)h(\vartheta)}{1+\vartheta^{\frac \alpha 2 - \frac 7{12}}}$, it follows that $f_\xi\rightarrow f$ almost everywhere. Further, using interpolation between \eqref{p142} and \eqref{p146} we get
$$
\|\vartheta_\xi\|_{L^{\alpha + \frac 32}((0,T)\times B)} \leq c
$$
  and $f_\xi\in L^{2 + \frac 1{\frac 32\alpha + \frac 74}}((0,T)\times B)$ uniformly. Consequently, $|f_\xi|^2$ is uniformly integrable and from the Vitali convergence theorem we get \eqref{convergence1}.\\
Finally, we have
\begin{equation*}
\chi_{\nu,\xi}(t,x) \kappa(\vartheta_\xi) h(\vartheta_\xi)\nabla \vartheta_\xi = \frac{\chi_{\nu,\xi}(t,x) \kappa(\vartheta_\xi) h(\vartheta_\xi)}{1+\vartheta_\xi^{\frac {\alpha}{2} - \frac 7{12}}}\left(\nabla \vartheta_\xi + \frac{1}{\frac \alpha 2 + \frac 5{12}}\nabla (\vartheta_\xi^{\frac \alpha 2 + \frac 5{12}})\right) =: a_\xi b_\xi.
\end{equation*}
Since $a_\xi$ tends strongly to $\frac{\chi_{\nu}(t,x) \kappa(\vartheta) h(\vartheta)}{1+\vartheta^{\frac {\alpha}{2} - \frac 7{12}}}$ in $L^2((0,T)\times B)$ and $b_\xi$ tends to $\left(\nabla \vartheta + \frac{1}{\frac \alpha 2 + \frac 5{12}}\nabla (\vartheta^{\frac \alpha 2 + \frac 5{12}})\right)$ weakly in $L^2((0,T)\times B)$, we get
$$
\chi_{\nu,\xi}(t,x) \kappa(\vartheta_\xi) h(\vartheta_\xi)\nabla \vartheta_\xi \rightarrow \chi_{\nu}(t,x) \kappa(\vartheta) h(\vartheta)\nabla \vartheta\ \mbox{ in }\ \mathcal D'((0,T)\times B).
$$
We conclude, that for $\varphi$ such that $\nabla \varphi \cdot {\bf n}|_{\Gamma_t} = 0$ we get
\begin{multline*}
I_\xi\rightarrow -\int_0^T\int_B \chi_{\nu} \nabla\K_h(\vt)\nabla \varphi \, \dxdt = -\int_0^T\int_{\Omega_t} \chi_\nu \nabla \K_h(\vt)\nabla \varphi \, \dxdt - \int_0^T \int_{B\setminus \Omega_t}\chi_\nu \nabla \K_h(\vt)\nabla \varphi \, \dxdt\\=  \int_0^T\int_{\Omega_t} \chi_\nu  \K_h(\vt)\Delta \varphi \, \dxdt + \int_0^T \int_{B\setminus \Omega_t}\chi_\nu \K_h(\vt)\Delta \varphi \, \dxdt = \int_0^T\int_B \chi_{\nu} \K_h(\vt)\Delta \varphi \, \dxdt\\= \int_0^T \int_B \K_{\nu, h}(\vt)\Delta \varphi \, \dxdt.
\end{multline*}

It remains to handle the term
$$
\int_0^T\int_B \chi_{\nu,\xi} \kappa(\vt_\xi)h'(\vt_\xi)|\nabla_x \vt_\xi|^2 \, \dxdt.
$$
However, since $\chi_{\nu,\xi}\sil \chi_\nu$ a.e in $B$ as $\xi\sil 0$, we may proceed similarly as in \eqref{s11} and \eqref{s12}.

Hence, after passing with $\xi \sil 0$ the thermal energy inequality has the form
	\bFormula{p42lim}
	\int_0^T\int_{B} (\vr+\delta) Q_h(\vt) \partial_t \varphi + \vr Q_h(\vt) \vu \cdot \nabla \varphi +
	 \K_{\nu,h}(t,x,\vt) \Delta \varphi - \delta h(\vt)\vt^{\alpha + 1}\varphi \, \dxdt
	 \eF
	 \[	\leq \int_0^T\int_{B} (h(\vt)\vt p_\vt \Div \vu - h(\vt)(1-\delta)\tn{S}_\omega : \Grad \vu + h'(\vt)\kappa_{\nu}(t,x,\vt)|\nabla_x\vt|^2 ) \varphi \, \dxdt + \int_B(\vr_{0,\delta} + \delta) Q(\vt_{0,\delta}) \varphi(0) \, \dx
	\]
for any $\varphi \in C^\infty([0,T] \times \mathbb{R}^3),$ $\varphi \geq 0,$ $\varphi(T)=0$, $\nabla \varphi \cdot \vc{n}|_{\partial B}=0$, $\nabla \varphi \cdot \vc{n}|_{\Gamma_\tau}=0$ for any $\tau \in [0,T]$ and any $h \in C^\infty([0,\infty))$ such that \eqref{eq:h} holds.

The equation of continuity \eqref{s2}, momentum equation \eqref{s20} and the global total energy inequality \eqref{s43} remain in the same forms.

\subsection{Vanishing viscosity limit}
In this section we let $\omega\sil 0$ in order to get rid of the last integral in \eqref{s20}. Let $\{ \vr_\omega, \vu_\omega, \vt_\omega\}_{\omega> 0}$ be the solution constructed in the previous section. Let us recall that  the viscosity coefficient has the following form  (see \eqref{p2}, \eqref{p2_2})
\[
\mu_\omega = \left\{ \begin{array}{l} \mu = {\rm const} > 0 \ \mbox{in} \ Q_T , \\ \\
\mu_\omega \to 0 \ \mbox{a.e. in} \ ((0,T) \times B) \setminus Q_T. \end{array} \right.
\]
From \eqref{p11} we deduce that
\bFormula{s22}
\int_0^T \int_{\Omega_t} \left| \Grad \vu_\omega + \Grad^t \vu_\omega - \frac{2}{3} \Div \vu_\omega \tn{I} \right|^2 \, \dxdt \leq c
\eF
and
\[
\int_0^T \int_{ B \setminus \Omega_t } \mu_\omega \left| \Grad \vu_\omega + \Grad^t \vu_\omega - \frac{2}{3} \Div \vu_\omega \tn{I} \right|^2 \, \dxdt \leq c.
\]
The last estimate yields
\[
\int_0^\tau \int_{ B \setminus \Omega_t } \mu_\omega \left( \Grad \vu_\omega + \Grad^t \vu_\omega - \frac{2}{3} \Div \vu_\omega \tn{I} \right) : \Grad \vph  \, \dxdt =
\]
\[
\int_0^\tau \int_{ B \setminus \Omega_t } \sqrt{\mu_\omega} \sqrt{\mu_\omega} \left( \Grad \vu_\omega + \Grad^t \vu_\omega - \frac{2}{3} \Div \vu_\omega \tn{I} \right)
: \Grad \vph \, \dxdt \to 0 \ \mbox{as}\ \omega \to 0
\]
for any fixed $\vph$.

As we know from Lemma \ref{La1}, the density $\varrho_\omega$ is supported by the ``fluid'' region $Q_T$. We use \eqref{p8}, \eqref{s22} and Korn's inequality to infer
\[
\int_0^T \int_{\Omega_t} | \Grad \vu_\omega |^2 \, \dxdt \leq c.
\]

By the arguments of Section \ref{s:penlimit}, we let $\omega \to 0$ and obtain the same form of the continuity equation \eqref{s2}. The momentum equation takes the form
\bFormula{s23}
\int_{\Omega_\tau} \vr \vu \cdot \vph (\tau, \cdot) \, \dx -  \int_{\Omega_0} (\vr \vu)_{0, \delta} \vph (0, \cdot) \, \dx
\eF
\[
= \int_0^\tau \int_{\Omega_t } \Big( \vr \vu \cdot \partial_t \vph + \vr [\vu \otimes \vu] : \Grad \vph  + p(\vr,\vt) \Div \vph
+ \delta \vr^\beta \Div \vph - \tn{S}(\Grad \vu) : \Grad \vph \Big) \, \dxdt
\]
for  any test function $\vph$ as in (\ref{s4a}). We would like to emphasise that the compactness of the density is necessary only in the ``fluid'' part $Q_T$ so a possible loss of regularity of $\vu_\omega$ outside $Q_T$ is irrelevant.

The thermal energy equation \eqref{p42lim} can be written in a form
\begin{multline}
	\int_0^T\int_{\Omega_t} (\vr_\omega+\delta) Q_h(\vt_\omega) \partial_t \varphi  +
	 \K_{\nu,h}(\vt_\omega) \Delta \varphi - \delta h(\vt_\omega)\vt_\omega^{\alpha + 1}\varphi  +  \vr_\omega Q_h(\vt_\omega)\vu_\omega \cdot \nabla \varphi \, \dx \dt \\ +
   \int_0^T\int_{B\setminus\Omega_t} \delta Q_h(\vt_\omega) \partial_t \varphi  +
	 \K_{\nu,h}(\vt_\omega) \Delta \varphi - \delta h(\vt_\omega)\vt_\omega^{\alpha + 1}\varphi \, \dx\dt
		\\\leq  \int_0^T\int_{\Omega_t}\left(h(\vt_\omega)\vt_\omega p_\vt \Div \vu_\omega  +  h(\vt_\omega)(\delta - 1)\tn{S}_\omega : \Grad \vu_\omega + h'(\vt_\omega)\kappa_\nu(\vt_\omega)|\nabla_x\vt_\omega|^2 \right) \varphi  \, \dx\dt \\   + \int_0^T\int_{B\setminus \Omega_t} \left( h(\vt_\omega)(\delta - 1)\tn{S}_\omega : \Grad \vu_\omega + h'(\vt_\omega)\kappa_\nu(\vt_\omega)|\nabla_x\vt_\omega|^2 \right) \varphi \, \dxdt + \int_B(\vr_{0,\delta} + \delta) Q(\vt_{0,\delta}) \varphi(0) \, \dx	\end{multline}
		for any $\varphi \in C^\infty([0,T] \times \mathbb{R}^3),$ $\varphi \geq 0,$ $\varphi(T)=0$, $\nabla \varphi \cdot \vc{n}|_{\partial B}=0$, $\nabla \varphi \cdot \vc{n}|_{\Gamma_\tau}=0$ for any $\tau \in [0,T]$ and any $h \in C^\infty([0,\infty))$ such that \eqref{eq:h} holds.
		Note, that the last integral on the right hand side is negative and thus can be omitted. 
		
		By the same arguments as in Section \ref{s:penlimit} 
		we get after passing with $\omega$ to zero
\begin{multline}\label{thermal.energy.inequality}
	\int_0^T\int_{\Omega_t} (\vr+\delta) Q_h(\vt) \partial_t \varphi  +
	 \K_{\nu,h}(\vt) \Delta \varphi - \delta h(\vt)\vt^{\alpha + 1}\varphi  +  \vr Q_h(\vt)\vu \cdot \nabla \varphi \, \dx \dt \\ + \int_0^T\int_{B\setminus\Omega_t} \delta Q_h(\vt) \partial_t \varphi  +\K_{\nu,h}(\vt)\Delta\varphi - \delta h(\vt)\vt^{\alpha + 1}\varphi \, \dx\dt
		\\\leq  \int_0^T\int_{\Omega_t}\left(h(\vt)\vt p_\vt \Div \vu  +  h(\vt)(\delta - 1)\tn{S} : \Grad \vu + h'(\vt)\kappa_\nu(\vt)|\nabla_x\vt|^2 \right) \varphi  \, \dxdt + \int_B (\vr_{0,\delta} + \delta) Q(\vt_{0,\delta}) \varphi(0) \, \dx .
\end{multline}

The total energy inequality \eqref{s43} can be written in the following way
	\begin{multline}	\int_{\Omega_{\tau}} \left( \vr_\omega\frac{1}{2} |\vu_\omega|^2 + \vr_\omega P_e(\vr_\omega) +
	\frac{\delta}{\beta-1}\vr_\omega^\beta + (\vr_\omega + \delta)Q(\vt_\omega) \right)(\tau,\cdot) \,\dx
	+ \int_{B\setminus \Omega_\tau} \delta Q(\vt_\omega(\tau,\cdot)) \, \dx \\
	+ \delta\int_0^\tau \int_{\Omega_t} \tn{S}_\omega:\Grad \vu_\omega + \vt_\omega^{\alpha+1} \,\dx \dt
	+ \int_0^\tau \int_{B\setminus\Omega_t} \delta \left(\tn{S}_\omega : \Grad \vu_\omega
	+ \vt_\omega^{\alpha + 1}\right) \, \dx\dt \\
	\leq \int_{\Omega_0}  \frac{1}{2} \frac{|(\vr\vu)_{0,\delta}|^2}{\vr_{0,\delta}} + \vr_{0,\delta} P_e(\vr_{0,\delta})
	+ \frac{\delta}{\beta-1}\vr_{0,\delta}^\beta + (\vr_{0,\delta} + \delta)Q(\vt_{0,\delta}) \,\dx
	+ \int_{\Omega_\tau}  (\vr_\omega \vu_\omega \cdot \vc{V}) (\tau, \cdot)\, \dx
	- \int_{\Omega_0}(\vr \vu)_{0,\delta} \cdot \vc{V}(0, \cdot) \, \dx\\
	+ \int_0^\tau \int_{\Omega_t}
	\Bigg(
	 \mu_\omega \left(\Grad \vu_\omega
	+ \Grad^t \vu_\omega - \frac{2}{3} \Div \vu_\omega \tn{I} \right) : \Grad \vc{V}
	- \vr_\omega \vu_\omega \cdot \partial_t \vc{V} - \vr_\omega \vu_\omega \otimes \vu_\omega : \Grad \vc{V}
	- p(\vr_\omega,\vt_\omega) \Div \vc{V} \\
	- \frac{ \delta }{\beta - 1} \vr_\omega^\beta \Div \vc{V}
	\Bigg)  \, \dx\dt
  	+ \int_{B\setminus\Omega_0} \delta Q(\vt_{0,\delta}) \,\dx
	+ \int_0^\tau \int_{B\setminus \Omega_t} \mu_\omega \left(\Grad \vu_\omega + \Grad^t \vu_\omega
	- \frac 23 \Div \vu_\omega \tn{I}\right): \Grad \vc{V} \, \dx \dt.
\end{multline}

On the right hand side we pass to the limit in the same way as in the momentum equation. The term $\int_0^\tau\int_{B\setminus \Omega_t}\delta \tn{S}_\omega : \Grad \vu_\omega$ on the left hand side can be neglected since it is positive and the rest of the terms converge due to the same arguments as in Section~\ref{s:penlimit}. This way as $\omega \to 0$ we obtain
\begin{multline} \label{energy.inequality.without.omega}	
\int_{\Omega_{\tau}} \left( \vr\frac{1}{2} |\vu|^2 + \vr P_e(\vr) + \frac{\delta}{\beta-1}\vr^\beta
+ (\vr + \delta)Q(\vt) \right)(\tau) \, \dx + \int_{B\setminus\Omega_{\tau}}\delta Q(\vt(\tau,\cdot))\, \dx \\+ \delta\int_0^\tau \int_{\Omega_t} \tn{S}:\Grad \vu + \vt^{\alpha+1} \, \dx \dt
+ \delta\int_0^\tau \int_{B\setminus\Omega_t} \vt^{\alpha + 1}\,  \dx \dt
\\ \leq \int_{\Omega_0}  \frac{1}{2} \frac{|(\vr\vu)_{0,\delta}|^2}{\vr_{0,\delta}} + \vr_{0,\delta} P_e(\vr_{0,\delta}) + \frac{\delta}{\beta-1}\vr_{0,\delta}^\beta + (\vr_{0,\delta} + \delta)Q(\vt_{0,\delta}) \, \dx
+ \int_{\Omega_\tau}   (\vr \vu \cdot \vc{V}) (\tau, \cdot) \ \dx - \int_{\Omega_0}(\vr \vu)_{0,\delta} \cdot \vc{V}(0, \cdot)
\, \dx\\
+ \int_0^\tau \int_{\Omega_t} \left( \mu \left(\Grad \vu + \Grad^t \vu - \frac{2}{3} \Div \vu \tn{I} \right) : \Grad \vc{V} - \vr \vu \cdot \partial_t \vc{V} - \vr \vu \otimes \vu : \Grad \vc{V} -
p(\vr,\vt) \Div \vc{V} - \frac{ \delta }{\beta - 1} \vr^\beta \Div \vc{V} \right)  \, \dx \dt
\\+ \int_{B\setminus\Omega_0} \delta Q(\vt_{0,\delta}) \, \dx.
\end{multline}

\subsection{Vanishing artificial pressure}

In this section we make the final limit procedure $\delta \sil 0$ while setting $\nu = \nu(\delta)$ in a suitable way. In order to derive estimates independent of $\delta$, which we need further, we proceed in the following way. Starting from \eqref{energy.inequality.without.omega} we get
\begin{multline}
\int_{\Omega_\tau} \left( \frac{1}{2} \vr |\vu|^2 + \vr P_e (\vr) + \frac{\delta}{\beta - 1} \vr^\beta + (\vr + \delta) Q(\vt) \right)(\tau, \cdot) \, \dx\\
\leq c(\vc{V})\left(1+ \| \Grad \vu\|_{L^2(Q_T)}
 +\int_0^\tau\int_{\Omega_t}\frac 12 \vr |\vu|^2 + \vr P_e(\vr) + \vr Q(\vt) \, \dx \dt\right),
\end{multline}
where $c(\vc{V})$ may depend also on $B$, $T$ and initial conditions but is independent of $\delta$. Using the Gronwall inequality, we get
\begin{multline}\label{odhad.na.zbytek}
{\rm ess}\sup_{\tau\in (0,T)}\left(\int_{\Omega_t} \left( \frac{1}{2} \vr |\vu|^2 + \vr P_e (\vr) + \frac{\delta}{\beta - 1} \vr^\beta + (\vr + \delta) Q(\vt) \right)(\tau, \cdot) \, \dx   + \int_{B\setminus\Omega_\tau}\delta Q(\vt)(\tau,\cdot)\, \dx\right)\\+ \delta \int_0^T \int_{B} \vt^{\alpha + 1} \,\dxdt\leq c \left(1+ \|\Grad \vu\|_{L^2(Q_T)}\right) .
\end{multline}
In the same way as in \cite{EF71}  one may derive  from equation \eqref{thermal.energy.inequality} that
\begin{multline} \label{odhad.na.gradient.rychlosti}
\|\Grad \vu\|^2_{L^2(Q_T)}\leq\int_0^T\int_{\Omega_t} \tn{S}:\Grad \vu \, \dx\dt \leq \int_0^T\int_{\Omega_t} \vt p_{\vt}(\vr)|\Div \vu| \, \dx \dt + c(1+\|\Grad \vu\|_{L^2(Q_T)})\\
\leq  \left(\int_0^T\int_{\Omega_t} \vt^2 p^2_{\vt}(\vr) \, \dx \dt\right)^{\frac 12}\|\Grad \vu\|_{L^2(Q_T)} + c(1+\| \Grad \vu\|_{L^2(Q_T)}).
\end{multline}
Further,
\begin{multline} \label{2.odhad.na.gradient}
\int_0^T\int_{\Omega_t} \vt^2 \vr^{\frac{2\gamma}3} \, \dx \dt
\leq \int_0^T\left(\int_{\Omega_t} \vr\vt^{\frac \alpha 2} \, \dx \right)^{\frac 4\alpha}
 \left(\int_{\Omega_t} \vr^{\frac {2\alpha \gamma - 12}{3\alpha - 12}} \, \dx\right)^{\frac {\alpha - 4} \alpha} \, \dt \\
\leq c\left({\rm ess}\sup_{t\in (0,T)} \int_{\Omega_t} \vr Q(\vt) \, \dx\right)^{\frac 2\alpha} \int_0^T \left(\int_{\Omega_t} \vr^\gamma \ \dx\right)^{\frac {2\gamma \alpha - 12}{3\gamma \alpha}} \, \dt \\
\leq c(1+ \| \Grad \vu\|_{L^2((0,T)\times {\Omega_t})})^{
\frac 4\alpha + \frac{2\gamma\alpha - 12}{3\gamma \alpha}}
\end{multline}
provided $\alpha \geq \frac{12(\gamma - 1)} \gamma$. Note that $\frac 4\alpha + \frac{2\gamma\alpha - 12}{3\gamma \alpha} < 2$ for any $\alpha\geq4$ and $\gamma\geq \frac 32 $. Combining this and \eqref{odhad.na.gradient.rychlosti} we obtain
\begin{equation}
\|\Grad \vu\|_{L^2(Q_T)} \leq c
\end{equation}
with constant independent of $\delta$. Going back to \eqref{odhad.na.zbytek} we get

\bFormula{p8delta}
{\rm ess} \sup_{t \in (0,T)} \| \sqrt{\vr} \vu (t , \cdot)  \|_{L^2 (Q_T)} \leq c,
\eF
\bFormula{p9delta}
\rm {ess} \sup_{t \in (0,T)} \left\| \vr (t, \cdot) \right\|_{L^\gamma(\Omega_t)} \leq c,
\eF
\bFormula{p10delta}
{\rm ess} \sup_{t \in (0,T)} \delta \| \vr (t, \cdot) \|^\beta_{L^\beta(\Omega_t)} \leq c,
\eF
\begin{equation}\label{p11delta}
\delta \int_0^T \int_{B} \vt^{\alpha+1} \, \dx \dt\leq c,
\end{equation}
\begin{equation}
\delta\int_0^T\int_{B}  Q(\vt) \, \dx \dt \leq c,
\end{equation}
with $c$ independent of $\delta$.

Let $\{ \vr_\delta, \vu_\delta, \vt_\delta \}_{\delta >0}$ be a solution constructed in the previous section as a consequence of
letting $\omega \to 0.$
Now we pass to the limit with $\delta \to 0$ in the weak formulations of the equations and inequalities \eqref{s2}, \eqref{s23}, \eqref{thermal.energy.inequality}, \eqref{energy.inequality.without.omega}. The key point here is again the strong convergence of densities which is obtained using the results on propagation of oscillations mentioned already in Section \ref{s:scd}, for more details see \cite{EF70}. Having this, we pass to the limit in the continuity equation \eqref{s2} which together with Lemma \ref{La1} yields \eqref{m1}. We also easily pass to the limit in the momentum equation \eqref{s23} with test function $\vph$ as in (\ref{s4a}). Concerning the total energy inequality \eqref{energy.inequality.without.omega} we proceed in the same steps as in \cite[Section 7.5]{EF70} to obtain \eqref{i3a3}. It remains to pass to the limit in the thermal energy inequality.


\subsubsection{Thermal energy inequality}
In order to get a proper limit in thermal energy equation, we set $\nu = \delta^2$, we take $h(\vt_\delta) = \frac 1{(1+\vt_\delta)^z}$, $z\in (0,1)$ in \eqref{thermal.energy.inequality} and we let $\delta$ to $0$. In order to prove that the terms integrated over the solid part $B \setminus \Omega_t$ vanish in the limit we proceed in the following way. We have, due to \eqref{p11delta} and the  H\"older inequality,
\begin{multline*}
\int_0^T\int_{B} \delta h(\vt_\delta) \vt_\delta^{\alpha + 1} \varphi \, \dx \dt
\leq c \int_0^T\int_{B}\delta \vt_\delta^{\alpha + 1-z} \, \dx \dt
= c \int_0^T\int_{B}\delta^{\frac{\alpha + 1 - z}{\alpha + 1}} \vt_\delta^{\alpha + 1-z} \delta^{\frac{z}{\alpha + 1}} \, \dx \dt\\
\leq c\delta^{\frac z{\alpha + 1}}  \left(\int_0^T\int_{B}\delta \vt_\delta^{\alpha+1}\, \dx \dt\right)^{\frac{\alpha + 1-z}{\alpha + 1}} \leq c\delta^{\frac z{\alpha + 1}} \rightarrow  0 \quad \mbox{ as } \delta \to 0.
\end{multline*}
Similarly,
\begin{multline*}
\left|\int_0^T\int_{B} \delta Q_h(\vt_\delta) \partial_t \varphi \, \dx \dt\right| \leq c \left(\int_0^T\int_{B}\delta \vt_\delta \, \dx \dt
+ \int_0^T\int_{B}\delta \vt_\delta^{\frac \alpha 2}\, \dx \dt\right) \\
\leq c \left(\delta^{\frac{\alpha}{\alpha + 1}} \left(\int_0^T\int_{B }\delta\vt_\delta^{\alpha + 1}\, \dx \dt\right)^{\frac 1{\alpha + 1}} + \delta^{\frac 12}\left(\int_0^T\int_{B}\delta \vt_\delta^{\alpha+1}\ \dx \dt\right)^{\frac 12}\right) \rightarrow 0 \quad \mbox{ as } \delta \to 0 .
\end{multline*}
Finally, by the choice $\nu = \delta^2$ we have
\begin{equation*}
\begin{split}
\left|\int_0^T\int_{B\setminus \Omega_t} \K_{\nu,h}(\vt_\delta)\Delta \varphi    \ ,\dx\dt \right|
& \leq \int_0^T\int_{B\setminus \Omega_t} \delta^2 \int_0^{\vt_\delta} \kappa(z) \ {\rm d}z |\Delta\varphi| \, \dx\dt
\\
&
\leq c\int_0^T\int_{B\setminus \Omega_t} \delta^2 \left(\vt + \vt^{\alpha +1}\right) |\Delta \varphi| \, \dx\dt  \to 0 \quad \mbox{ as } \delta \to 0 .
\end{split}
\end{equation*}
We omit the details concerning the limit in the terms integrated on the fluid part $\Omega_t$, since it is done in the same way as in \cite[Section 7.5.5]{EF70} by first proceeding with $\delta \sil 0$ keeping $z$ fixed and then passing with $z \sil 0$. In the limit, \eqref{m6} is obtained, which completes the proof of Theorem \ref{Tm1}.

\section{Discussion}
\label{d}

The assumption on monotonicity of the pressure is not necessary, the same result can be obtained for a non-monotone pressure adopting the method
developed in \cite{EF61}.


The case of nonconstant viscosities $\mu = \mu(\vt)$ and $\eta = \eta(\vt)$ for bounded, continuously differentiable, globally Lipschitz functions $\mu,\eta$ with $\mu$ bounded away from zero, can be also treated by adopting the technique of \cite{EF71}.

As pointed out in the introduction, the general Navier slip conditions (\ref{b1}) are obtained introducing another boundary integral in the weak formulation, namely
\[
\int_0^T \int_{\Gamma_t}  \zeta (\vu - \vc{V}) \cdot \vph \ {\rm dS}_x \ \dt
\]
Taking $\zeta = \zeta(x)$ as a singular parameter, we can deduce results for mixed type no-slip - (partial) slip
boundary conditions prescribed on various components of $\Gamma_t$.


\def\ocirc#1{\ifmmode\setbox0=\hbox{$#1$}\dimen0=\ht0 \advance\dimen0
  by1pt\rlap{\hbox to\wd0{\hss\raise\dimen0
  \hbox{\hskip.2em$\scriptscriptstyle\circ$}\hss}}#1\else {\accent"17 #1}\fi}

\end{document}